\documentclass{amsart}
\usepackage{graphicx}
\newtheorem{thm}{Theorem}[section]
\newtheorem{cor}[thm]{Corollary}
\newtheorem{lem}[thm]{Lemma}

\theoremstyle{definition}

\theoremstyle{remark}

\numberwithin{equation}{section}

\begin{document}

\title[Metrics of constant scalar curvature]{Metrics of
constant scalar curvatures conformal to a  Riemannian
product with a round sphere}%
\author{Jimmy Petean}%
\address{Calle de Abajo 12, San Javier, Guanajuato, Mexico}%
\email{jimmy@cimat.mx}%

\thanks{me}%
\subjclass{53}%
\keywords{scalar curvature}%

\begin{abstract}
We consider the conformal class of the Riemannian product $g_0 + g$,
where $g_0$ is the constant curvature metric on $S^m$ and $g$ is a
metric of constant scalar curvature on some closed manifold. We show
that the number of metrics of constant scalar curvature in the
conformal class grows at least linearly with respect to the square
root of the scalar curvature of $g$. This is obtained by studying
radial solutions of the equation $\Delta u -\lambda u + \lambda u^p
=0$ on $S^m$, and the number of solutions in terms of $\lambda$.
\end{abstract}
\maketitle
\section{Introduction}

Any closed  manifold admits metrics of constant scalar curvature.
Given any Riemannian metric $g$ on $M^n$ we consider its conformal
class $[g]$, and define the Yamabe constant of $[g]$ as the minimum
of the (normalized) total scalar curvature functional restricted to
$[g]$,

$$Y(M,[g]) = \inf_{h\in [g]} \frac{\int_M s_h dvol_h
}{Vol(M,h)^{\frac{n-2}{n}}} $$

\noindent where $s_h$ and $dvol_h$ are the scalar curvature and
volume element of $h$.

It is elementary that the functional restricted to $[g]$ is bounded
below and the fact that the infimum is actually achieved is a
fundamental result, obtained in a series of steps by Hidehiko Yamabe
\cite{Yamabe}, Thierry Aubin \cite{Aubin}, Neil Trudinger
\cite{Trudinger} and Richard Schoen \cite{S2}. Since the critical
points of the functional (restricted to $[g]$) are the metrics of
constant scalar curvature in $[g]$ it follows that minimizers are
metrics of constant scalar curvature. These are called {\it Yamabe
metrics}. So in any conformal class of metrics in any closed
manifold there is at least one unit volume metric of constant scalar
curvature. If the Yamabe constant of $[g]$ is non-positive there is
actually only one, the Yamabe metric of the conformal class. But
when the Yamabe constant is positive there might be more. For
instance Daniel Pollack proved in \cite{Pollack} that every
conformal class with positive Yamabe constant can be $C^0$
approximated by a conformal class with an arbitrarily large number
of (non-isometric) metrics of constant scalar curvature. But to
determine all the metrics of constant scalar curvature in a given
conformal class of positive Yamabe constant is a very difficult
problem. Uniqueness still holds for the conformal class of positive
Einstein metrics different from the round metric on $S^n$ by a
result of Morio Obata \cite{Obata}. For the conformal class of the
round metric all constant scalar curvature metrics in the conformal
class are obtained by conformal diffeomorphisms of the sphere (a
non-compact family) and are all isometric. Examples of multiplicity
of metrics of constant scalar curvature in a conformal class (which
is a Riemannian covering of a number of manifolds) are also obtained
by Emmanuel Hebey and Michel Vaugon in \cite{Hebey}. Recently, Simon
Brendle gave examples of smooth conformal classes of Riemannian
metrics on high dimensional spheres for which the space of unit
volume constant scalar curvature metrics in the conformal class is
non-compact \cite{Brendle}.

For any closed Riemannian manifold $(M^n ,g)$ the Yamabe constant of
its conformal class is bounded above by $Y(S^n ,[g_0 ])$, where
$g_0$ is the round metric \cite{Aubin}. Therefore if $(M_1 , g_1 )$,
$(M_2 ,g_2 )$ are Riemannian manifolds of constant scalar curvature
and $s_{g_1 }$ is positive then for $\delta$ positive and small the
conformal class of the Riemannian product $\delta g_1 + g_2$ has at
least two constant scalar curvature metrics: the Riemannian product
and a Yamabe metric. The simplest case is to consider the Riemannian
product $g_0 + dt^2$ on $S^n \times S^1$. In this case (see the
article by R. Schoen \cite{S3} and Osamu Kobayashi \cite{Kobayashi1,
Kobayashi2}) all conformal factors producing metrics of constant
scalar curvature are functions of $S^1$, and there is a sequence of
values $\delta_i \rightarrow 0$ such that for $\delta \in
(\delta_{i+1} , \delta_i )$ the number of constant scalar curvature
metrics in the conformal class of $\delta g_0 +dt^2$ is $i$.

In this article we will draw a similar picture in the case of $S^k
\times S^m$, $k,m>1$. As usual we will call $p=p_N = \frac{2N}{N-2}$
and $a=a_N = \frac{4(N-1)}{N-2}$. We will consider the Riemannian
product $\delta g_0^k +g_0^m$ which has scalar curvature ${\bf
s}_{\delta} = (1/\delta ) (k(k-1)) + m(m-1)$. For a positive, radial
function $f:S^m \rightarrow {\bf R}$, let $u:[0,\pi ] \rightarrow
{\bf R}$ be the corresponding function (so $f(x) = u(d(x,P))$, where
$P$ is a fixed point in $S^m$). The conformal metric
$f^{\frac{4}{k+m-2}} \ (\delta g_0^m + \delta_0^k )$ has constant
scalar curvature $K$ if and only if $f$ satisfies the Yamabe
equation

$$-a_{m+k} \Delta_{g_0^m} \ f+ {\bf s}_{\delta} f = K f^{p_{m+k} -1}.$$

\noindent Then the function $u$ must satisfy the equation

$$u''+(m-1) \frac{\cos (t)}{\sin (t)} u' + \frac{K}{a} u^{p-1}
-\frac{{\bf s}_{\delta}}{a} u =0.$$

We normalize taking $K={\bf s}_{\delta}$, and call $\lambda = K/a$.
Then we are looking for positive solutions of the equation

$$u'' + (m-1) \frac{\cos (t)}{\sin (t)} u' + \lambda (u^{p-1} -u) =0$$

\noindent with initial conditions $u(0)=\alpha$, $u'(0)=0$ and such
that $u'(\pi )=0$.

The corresponding equation in ${\bf R}^m$ has been well studied.
First one has to note that in this case the equations for different
values of $\lambda$ are all equivalent. So one only has to consider
the equation $\Delta f -f +f^q =0$. Then from the  classical work of
Basilis Gidas, Wei-Ming Ni and Louis Nirenberg \cite{Gidas, Ni}
follows that all solutions which are positive and vanish at $\infty$
(ground states) must be radially symmetric. Then one is looking for
solutions of the equation

$$u''+\frac{m-1}{t} u' + u^q -u=0,$$

\noindent with initial conditions $u(0)=\alpha$, $u'(0)=0$ which are
positive and $u(\infty )=0$. This equation has been completely
analyzed by Man Kam Kwong in \cite{Kwong}, proving in particular
that there exists exactly one such solution.

In our case we will see that the number of solutions grows at least
linearly in $\sqrt{\lambda}$. We will build solutions which verify
that $u'(\pi /2)=0$ so the corresponding metric is invariant under
the antipodal map, producing a metric in the projective space. This
will therefore prove the following:

\begin{thm} Let $(M^k ,g)$ be a Riemannian  manifold of constant scalar curvature ${\bf s}$.
The number of unit volume non-isometric metrics of constant scalar
in the conformal conformal class $[g_0 +g]$ on $S^m \times M$ grows
at least linearly with $\sqrt{\bf s}$. The same is true if we
replace $S^m$ with the projective space ${\bf P}^m$, with the metric
of constant curvature. More explicitly, we will show that if $n\geq
1$ and

$$({\bf s} +m(m-1))\ \left( \frac{p-2}{a} \right) \in ( 2n(2n+m-1), (2n+2)(2n+2+m-1) ]$$

\noindent then $[g_0 +g]$ contains at least $2n +2$ unit volume
non-isometric metrics of constant scalar curvature.
\end{thm}

In order to construct the solutions in the previous theorem we will
need to prove that there exists one radial solution which is
strictly decreasing in $[0,\pi ]$. To do so we will need the
following elementary result, which we will prove in Section 4:

\begin{thm} Let $(M^k ,g)$ be a closed Riemannian manifold of
constant scalar curvature ${\bf s}$. If

$${\bf s} +m(m-1) > \frac{a \ m}{p_{k+m}-2} $$

\noindent then the Riemannian product $g+g^m_0$ is not a Yamabe
metric. Actually, the product metric is not a local minimum of the
total scalar curvature functional restricted to the space $\{
f(g+g_0^m ): f:S^m \rightarrow {\bf R}_{>0} \}$.
\end{thm}

\vspace{.5cm}

Let us consider the case $m=k=2$. Then $p=p_4 = 4$ and $a=a_4 =6$.
We study the equation $u''+(\cos (t) / \sin (t)) u' + \lambda (u^3
-u)=0$ and $\lambda$ relates to the (constant) scalar curvature of
$(M,g)$ as $\lambda=(1/6)({\bf s} +2)$. Let $A_n = (1/2)n(n+1)$. We
will show that for $\lambda \in (A_1 , A_2 ]$ there at least two
solutions; one of them is the constant solution which is not a
Yamabe minimizer and the other one is a strictly decreasing
function. For $\lambda \in (A_2 , A_3 ]$ there are at least 4
solutions and in general for $\lambda \in (A_{2n} , A_{2n+2} ]$
there are at least $2n+2$ solutions. Except from one of them, all of
these solutions verify that $u' (\pi /2 )=0$ and so produce constant
scalar curvature metrics on ${\bf P}^2 \times M$.

Probably the most interesting particular case is to consider the
conformal classes of the Riemannian product of metrics of constant
curvature  on $S^2 \times S^2$, which we will write $g_0 + \delta
g_0$. Moving $\delta$ in $(0, \infty )$ we obtain the values of
$\lambda$ in the range $(1/3, \infty )$. The previous comments
translate into the following:

\begin{thm} The metric $g_0 +\delta g_0$ is not a Yamabe metric for
$ \delta <1/2$. Let $\delta_n = 2\ (3n(n+1)-2)^{-1}$. For $\delta
\in [\delta_{2(n+1)} , \delta_{2n} )$ the number of constant scalar
curvature metrics in the conformal class of $g_0 +\delta g_0$ is at
least $2n+2$.
\end{thm}

{\bf Remark}: For $\delta =1/2$ the conformal class of $g_0 + (1/2)
g_0$ attains the same value of the Yamabe functional, $12\sqrt{2}
\pi$, as the conformal class of the Fubini-Study metric on ${\bf
CP}^2$. It is hard to imagine that this is accidental and one would
expect that there is a geometric explanation for the coincidence.

{\bf Remark}: The results of Gidas, Ni, Nirenberg on the symmetry of
solutions in ${\bf R}^n$ do not seem to apply to the case of $S^n$,
and there is no adaptation of them in the literature, at least to
the author's knowledge. But it seems reasonable to expect that there
must be some variation of their arguments proving that solutions on
$S^n$ are all radially symmetric. All but one of the solutions we
are going to prove that exist have a local maximum at both 0 and
$\pi$ or a local minimum at both 0 and $\pi$. It is clear that in
between any two of these solutions there exists one solutions which
has a minimum at 0 and a maximum at $\pi$ (or viceversa). These
should be all the radially symmetric solutions, but to prove this
one should adapt many subtle ideas appearing  in the work of M. K.
Kwong \cite{Kwong}. This would describe all solutions which depend
on only one of the factors. All in all it seems reasonable to
conjecture that for $\delta \geq 1/2$ the metric $g_0 + \delta g_0$
is a Yamabe metric and the only unit volume metric of constant
scalar curvature in its conformal class (Obata's Theorem
\cite{Obata} says that this is true for $\delta =1$). And that for
$\delta <1/2 $ the previous comments describe all the constant
scalar curvature metrics in the conformal class for which the
conformal factor depends on only one of the spheres. It is tempting
to guess that these are actually all the solutions, but there is no
real evidence to support that.

\vspace{.3cm}

\noindent {\it Acknowledgements:} The author would like to thank
Claude LeBrun for very useful comments on the original draft of this
manuscript.

\section{Sturm comparison}

To study the differential equation we will need to apply some Sturm
comparison techniques. For the convenience of the reader we will
state the appropriate version of Sturm Theorem. It appears in
\cite[Lemma 1]{Kwong}.

\begin{thm} Let $U$ and $V$ be solutions of the equations

$$U''(t) +f(t) U'(t) + g(t)U(t) =0, \ \ t\in (a,b),$$

$$V''(t)+ f(t) V'(t) + G(t)V(t)=0, \ \ t\in (a,b).$$

\noindent Let $(\alpha ,\beta )$ be a subinterval where $V(t)\neq 0$
and $U(t)\neq 0$ and such that $G(t)\geq g(t)$ for all $t\in
(\alpha, \beta )$.

If

$$\frac{V'(\alpha )}{V(\alpha )} \leq  \frac{U'(\alpha )}{U(\alpha
)},$$

\noindent then

$$\frac{V'(t)}{V(t)} \leq \frac{U'(t)}{U(t)} \ \ \ \forall \ t\in
(\alpha ,\beta ).$$

\noindent If equality holds at any $x\in (\alpha, \beta )$ then
$U\equiv V$ in $[\alpha , x]$.

\end{thm}

\section{Solutions near $u=1$; the linear ODE}

We want to study the differential equation

$$u'' + (m-1) \frac{\cos (t)}{\sin (t)} u' + \lambda (u^{p-1} - u ) =0.$$

\noindent where $\lambda$ is positive, $p>2$ and $m-1$ is a positive
integer. We set the initial conditions to be $u(0)=\alpha $ and
$u'(0)=0$. The interval of definition is $[0,\pi ]$ and we are
interested in positive solutions such that $u' (\pi )=0$ (which give
solutions in $S^m$. We consider $u=u(t,\alpha , \lambda )$.

There is a canonical solution, $u(t,1,\lambda )=1$. Our goal in this
section is to understand the behavior of solutions near this
canonical one; solutions $u(t,\alpha ,\lambda )$ with $\alpha$ close
to 1.

Consider the function

$$w(t)=\frac{\partial u}{\partial \alpha} (t,1,\lambda ).$$

\noindent Then $w$ is solution of the linear equation

$$w'' + (m-1) \frac{\cos (t)}{\sin (t)} w' + \lambda (p-2) w =0 ,$$

\noindent with the initial conditions $w(0)=1, w'(0)=0$.

We let $A=(p-2) \lambda$ and call $w=w_A$ the corresponding
solution. The solutions for $A=n(n+m-1)$ can be given explicitly.

For instance

$$w_0 =1, \  w_m (t) = \cos (t), \  w_{2(m+1)} (t)=\frac{m+1}{m} \left( \cos^2 (t)
-\frac{1}{m+1} \right) .$$

 If we call

$$H_A (f)=f'' +(m-1)\frac{\cos (t)}{\sin (t)} f' + A f,$$

\noindent then we have

$$H_A (\cos^n (t))= (A- n(n+m-1))\cos^n (t)+ n(n-1) \cos^{n-2} (t).$$

It easily follows

\begin{lem} $w_{n(n+m-1)}$ is a linear combination (with rational
coefficients) of powers of $\cos^{n-2k} (t)$, where $0\leq 2k \leq
n$.
\end{lem}

Therefore it follows that if $n$ is odd then $w_{n(n+m-1)} (\pi )
=-1$, and if $n$ is even $w_{n(n+m-1)} (\pi ) = 1$. Moreover, if $n$
is even $w_{n(n+m-1)}$ is symmetric with respect to $t=\pi /2$ (and
therefore $w_{n(n+m-1)} (\pi /2 ) \neq 0$, since in that case by the
uniqueness of solutions it would have to vanish everywhere) and if
$n$ is odd $w_{n(n+m-1)}$ is antisymmetric with respect to $t=\pi
/2$ and (and $w_{n(n+m-1)} (\pi /2)=0$).

\begin{lem} The solution $w_{n(n+m-1)}$ has exactly $n$ zeros in the
interval $(0,\pi )$. The number of zeros in the interval $(0,\pi
/2)$ is equal to the number of zeros in the interval $(\pi /2,0)$.
\end{lem}

\begin{proof} We use induction  on $n$. We know it is true for the
first values of $n$ by explicit computation. By Sturm comparison
(Theorem 2.1) we know that if $A<B$ then the solution $w_B$ has at
least one 0 in between any two zeros of $w_A$. Therefore $w_B$ has
at least the same number of zeros as $w_A$, and if it has exactly
the same then both must have the same sign after the last 0. Since
when moving from n to n+1 the corresponding solutions change sign at
the final value $\pi$ it follows that $w_{(n+1)(m+n)}$ must have at
least one more 0 than $w_{n(n+m-1)}$. By induction this means than
$w_{(n+1)(m+n)}$ must have at least n+1 zeros. But $w_{(n+1)(m+n)}$
is a polynomial of degree n+1 in $\cos (t)$ and $\cos (t)$ is
injective in $(0,\pi )$. Therefore $w_{(n+1)(m+n)}$ could have at
most n+1 zeros.

The last statement follows directly from the previous comments.
\end{proof}

The information we will use to prove the existence of constant
scalar curvature metrics is about the number of local extrema of
$w_A$ in the interval $(0,\pi /2)$. We can give a complete analysis
of this. For $n$ even $w_{n(n+m-1)} '(\pi /2 )=0$, and the number of
local extrema in $(0,\pi /2)$ is $n/2 -1$ (from the previous lemma).
By Sturm comparison (Theorem 2.1) the number of local extrema of
$w_A$ is a non-decreasing function of $A$. This function jumps by
one every time we cross a value $A=2n(2n+m-1)$. Therefore if we call
$C_n = (2n(2n+m-1), (2n+2)((2n+2+m-1)]$, we have proved

\begin{thm} For $A\in C_n$, the solution $w_A$ has exactly $n$ local
extrema in $(0,\pi /2)$.
\end{thm}

\section{Proof of Theorem 1.2 and  the existence of a strictly decreasing
solution}

First recall that if we have two conformal metrics ,$H$ and $G$, on
an $N$-dimensional manifold and we express the conformal relation as
$H=f^{\frac{4}{N-2}} G$ then the expression for the total scalar
functional of $H$, $S(H)$, in terms of $G$ and $f$ is

$$S(H)=Y_G (f) =\frac{4a\int | \nabla f |^2 dvol_G \
+ \ \int s_G f^2 dvol_G }{ {\left( \int f^p dvol_G \right) }^{2/p}}
= \frac{E_G (f)}{{\| f \|}_p^2}.$$

\noindent Recall also that

$$(d/dt)|_{t=0} \  (Y_G (f+tu))=\frac{2}{{\| f\|}_p^2} \ \int [-a
\Delta f + sf - {\| f\|}_p^{-p} E_G (f) f^{p-1} ] u \ dvol_G .$$

\vspace{.5cm}

Given a Riemannian product of constant scalar curvature metrics $g_1
+g_2$ on $M_1 \times M_2$ one can consider conformal factors
depending on only one of the variables and define, for instance,
$[g_1 +g_2]_1 = \{ f. (g_1 + g_2 ) : f:M_1 \rightarrow {\bf R}_{>0}
\}$. Then we  define \cite{AFP}

$$Y_1 (M_1 \times M_2 ,g_1 + g_2 )= \inf_{h\in [g_1 + g_2 ]_1 } \
\frac{\int_{M_1 \times M_2} \ s_h \ dvol_h}{(Vol(M_1 \times M_2
,h)^{\frac{N-2}{N}}}, $$

\noindent where $N= dim(M_1 \times M_2 )$. It is easy to see that
the infimum is realized \cite[Proposition 2.2]{AFP}. In the case
$(M_1 ,g_1 )=(S^m ,g_0 )$ given any positive function $f$ on $S^n$
one can consider the spherical symmetrization $f_*$, which is the
radial non-increasing function on $S^n$ which verifies $Vol \{ f_*
>t \} = Vol \{ f>t\}$ for all $t>0$. Then it is well-known that the
total scalar curvature functional is non-increasing by this
symmetrization, namely $S(f_* (g_1 +g_2 )) \leq S(f (g_1 + g_2 )$.
This proves:

\begin{lem} If $(M^k ,g)$ has constant scalar curvature there exists
a radially symmetric non-increasing function on $S^m$ which gives a
minimizer for $Y_1 (S^m \times M ,g+ g_0 )$.
\end{lem}

\begin{proof} {\bf (Theorem 1.2)} As we mentioned in the previous
section, the function $u(x) =\cos (d(x, N))$ is an eigenfunction of
the (negative) Laplacian operator of $(S^m ,g_0 )$ (and hence of
$(M\times S^m ,g_0 +g)$ with eigenvalue $-m$. Moreover,

$$\int_{S^m } u \ dvol_{g_0} =0.$$

\noindent Let $Y(t)= Y_{g_0 +g} (1+tu)$. Then $Y' (0) =0$ and for
some positive constant $K$,

$$Y'' (0) = K \left( \int_{S^m} -a\Delta u \ u + (s_g + m(m-1)) u^2
-(p-1)(s_g +m(m-1)) u^2 \ dvol_{g_0} \right) $$

$$=K (a m +(2-p)(s_g +m(m-1))) \int_{S^m} u^2 \ dvol_{g_0}.$$

The hypothesis says precisely that the previous expression is
negative and this proves the theorem.

\end{proof}

\vspace{1cm}

Theorem 1.2 and Lemma 4.1 imply:

\begin{cor} If $\lambda >\frac{m}{p_{m+k}-2}$ then there is $\alpha >1$
such that the solution of the  equation

$$u''+(m-1)\frac{\cos (t)}{\sin (t)} u' +\lambda (u^{p-1} -u)=0$$

\noindent with initial conditions $u(0)=\alpha$, $u'(0) =0$, is
positive, strictly decreasing in $[0,\pi ]$, and $u'(\pi )=0$.
\end{cor}

\section{The number of solutions: proof of Theorem 1.1}

Now we fix $\lambda >0$  and study the dependence of the solution
$u(t,\alpha ,\lambda )$ of the equation

$$u'' + (m-1)\frac{\cos (t)}{\sin (t)} + \lambda (u^{p-1} -u) =0$$

\noindent  on $\alpha$. It follows from Section 1 that if $\alpha$
is close to 1, then $u$ intersects the canonical solution about
$\sqrt{\lambda }$ times.

Let $P$ be the subset $\{ \alpha \in (0,\infty ) : u_{\alpha }>0$ on
$[0,\pi /2] \}$. Clearly $P$ is an open subset of $(0,\infty )$ and
$1\in P$. Suppose $[1,A)$ is a maximal (to the right) interval
included in $P$. Then $u_A$ must be nonnegative on $[0,\pi /2]$.
Then $u_A$ must be strictly positive on $[0,\pi /2)$ and $u_A (\pi
/2) =0$ (otherwise the interval would not be maximal, of course).

Now consider the Energy function associated to $u_A$,

$$E_A (t) = \frac{(u_A'(t))^2}{2} + \lambda \left( \frac{u_A^{q+1} (t)}{q+1} -
\frac{u_A^2 (t)}{2} \right) .$$

We have

$$E_A' (t) = - m\frac{\cos (t)}{\sin (t)} (u_A' (t))^2.$$

And so $E_A$ is decreasing in the interval $[0, \pi /2]$. Since $E_A
(\pi /2) >0$ we must have positive energy on $[0,\pi /2]$. Consider
the following simple lemma:

\begin{lem} If $u_{\alpha}$ has a local minimum at $t_0$
then $u_{\alpha } (t_0 )<1$ and so $E_{\alpha} (t_0 ) <0$.
\end{lem}

Then it follows:

\begin{lem} If for some finite $A$, $[1,A)$ is a maximal (to the right)
interval contained in $P$ then $u_A$ has no local extrema in $(0,\pi
/2)$.
\end{lem}

The following lemma will allow us to construct solutions without
having to analyze the equation in the whole interval $[0,\pi ]$.

\begin{lem} Suppose that for some positive $\alpha $ the solution $u_{\alpha}$
verifies $\ u_{\alpha }' (\pi /2) =0$. Then $u_{\alpha}' (\pi )=0$
(and actually $u_{\alpha}$ is symmetric with respect to $t=\pi /2$).
\end{lem}

\begin{proof} The  function $v(t)= u_{\alpha} (\pi -t)$ for $t\in
[\pi /2 , \pi)$ is also a solution of the equation. Moreover $v(\pi
/2)= u_{\alpha} (\pi /2 )$ and $v' (\pi /2)=0=u_{\alpha}' (\pi /2
)$. Therefore $v=u_{\alpha}$ and the lemma follows.
\end{proof}

\begin{lem} Suppose that for some positive $\alpha_0 (\neq 1)$ the solution
$u_{\alpha_0}$ has exactly $k$ extrema in the open interval $(0,\pi
/2)$ and $u_{\alpha_0}' (\pi /2 )=0$. Then there exists
$\varepsilon>0$ such that for $\alpha \in (\alpha_0 -\varepsilon
,\alpha_0 +\varepsilon )$ the number of local extrema of
$u_{\alpha}$ in $(0, \pi /2 ]$ is $k$ or $k+1$.
\end{lem}

\begin{proof} If $u_0 =u_{\alpha_0} (\pi /2)=1$ we would have $\alpha_0 =1$, and
we have assumed this is not the case. Therefore $u_0 <1$ or $u_0
>1$. If $u_0 <1 (>1)$ there exists $\delta >0$ such that for $\alpha
\in (\alpha_0 -\delta ,\alpha_0 +\delta )$ and $t\in (\pi /2 -\delta
,\pi /2 + \delta )$, we have $u_{\alpha} (t) <1 (>1)$. This implies
that such $u_{\alpha}$ can not any local maximum (minimum) in $(\pi
/2 -\delta , \pi /2 +\delta )$. Therefore it has at most 1 local
extrema in that interval.

We can also assume the $\delta$ small enough so that $u_{\alpha_0 }$
does not have any other extrema besides $\pi /2$ in $[\pi /2 -\delta
, \pi /2 + \delta ]$. Therefore $u_{\alpha_0}$ has $k$ local extrema
in that interval and for $\varepsilon
>0$ small enough, $\varepsilon  <\delta$, and $\alpha \in (\alpha_0
-\varepsilon ,\alpha_0 + \varepsilon )$, $u_{\alpha}$ also has $k$
local extrema in $(0,\pi /2 -\delta )$ (and hence $k$ or $k+1$ in
$(0,\pi /2]$.

\end{proof}

\begin{lem} If $\alpha$ is close to zero the solution $u_{\alpha}$
has no local extrema in $(0,\pi /2)$. If $\lambda (p-2) >m$ there
exists $\alpha >1$ such that the solution $u_{\alpha}$ has no local
extrema in $(0,\pi )$.
\end{lem}

\begin{proof} For $\alpha$ close to 0 the solution $u_{\alpha}$ stays
close to 0 and so stays less than 1 until $\pi /2$ and so it does
not have any local maximum and is increasing. This proves the first
statement. The second statement is just Corollary 4.2.

\end{proof}

We are finally ready to prove Theorem 1.1

\begin{proof} If $\lambda (p-2) \in (2n(2n+m), (2n+2)(2n+2+m)]$, for $\alpha$ close
to 1, it follows from Theorem 3.3 that the solution $u_{\alpha}$ has
at least $n$ local extrema in $(0,\pi /2)$. Increasing $\alpha$ from
1 to $\infty$ we bump into solutions for which $u_{\alpha}' (\pi
/2)=0$. Each one of these gives a constant scalar curvature metric.
As we cross this value of $\alpha$ the number of local extrema
before $\pi /2$ decreases at most by 1, from Lemma 5.3. It follows
from Lemma 5.2 and Lemma 5.5 that eventually we reach a value of
$\alpha$ for which the corresponding solution does not have any
local extrema in $(0,\pi /2)$. We can make the same argument when
$\alpha$ decreases from 1 to 0. Therefore the number of initial
values $\alpha$ for which $u_{\alpha}' (\pi )=0$ is at least $2n$.
Note that for these solutions 0 and $\pi$ are both local minima or
maxima. Besides these, we have the constant solution and one
strictly decreasing (or increasing) solution.
\end{proof}

\end{document}